\documentclass[fleqn]{article}
\usepackage{amssymb}
\usepackage{amsmath}

\newtheorem{definition}{\bf Definition}[subsection]
\newtheorem{theorem}[definition]{Theorem}

\newtheorem{conclusion}[definition]{Conclusion}

\newtheorem{corollary}[definition]{Corollary}

\newtheorem{example}[definition]{Example}

\newtheorem{proposition}[definition]{Proposition}
\newtheorem{remark}[definition]{Remark}

\newenvironment{proof}[1][Proof]{\noindent\textbf{#1.} }{\ \rule{0.5em}{0.5em}}
\input{tcilatex}
\roman{enumii}

\numberwithin{equation}{section}
\begin{document}

\title{A Morita context and Galois extensions for Quasi-Hopf algebras}
\author{Adriana Balan \\
%EndAName
Department of Mathematics, University Politehnica of Bucharest, Romania}
\date{October, 2006}
\maketitle

\begin{abstract}
If $H$ is a finite dimensional quasi-Hopf algebra and $A$ is a left $H$%
-module algebra, we prove that there is a Morita context connecting the
smash product $A\#H$ and the subalgebra of invariants $A^{H}$. We define
also Galois extensions and prove the connection with this Morita context, as
in the Hopf case.
\end{abstract}

\section{Introduction}

The Hopf Galois extensions appear for the first time in the papers of Chase,
Harrison and Rosenberg (\cite{Chase65})\ and of Chase and Sweedler (\cite%
{Chase69}). The actual definition is due to Kreimer and Takeuchi (\cite%
{Kreimer81}). A first Morita context has been constructed by Cohen, Fischman
and Montgomery in \cite{Cohen90}. They start from a finite dimensional Hopf
algebra $H$ over a field and an $H$-module algebra $A$, and construct a
Morita context connecting the smash product $A\#H$ and the subalgebra of $H$%
-invariants $A^{H}$, showing also the connection with the Galois extensions.
Another Morita context has been constructed by Doi, generalizing a
construction of Chase and Sweedler (\cite{Chase69}). In this case, the
author began with a Hopf algebra $H$ not necessarily finite dimensional and
an $H$-comodule algebra $\mathcal{A}$. The rings connected by the Morita
context are this time $\mathcal{B}$, the subalgebra of coinvariants, and $%
\#(H,\mathcal{A})$.

In the case of a finite dimensional Hopf algebra $H$, a left $H$-module
algebra is the same as a right $H^{\ast }$-comodule algebra, so it is
natural to ask if both contexts coincide. The affirmative answer has been
given by Beattie, D\u{a}sc\u{a}lescu and Raianu in \cite{Beattie97a}, for a
co-Frobenius Hopf algebra.

The books of Montgomery \cite{Montgomery93} and of D\u{a}sc\u{a}lescu, N\u{a}%
st\u{a}sescu and Raianu \cite{Dascalescu00} are a good reference about the
main results on this subject.

On the other hand, in the last fifteen years, several different
generalizations of Hopf algebras have appeared: corings, weak Hopf algebras,
quasi-Hopf algebras and Hopf algebroids. The Galois theory for corings was
realized by Brz\`{e}zinski (\cite{Brzezinski02}), while the Morita theory
was developed for corings by Caenepeel, Vercruysse and Wang (\cite%
{Caenepeel04}). The case of weak Hopf algebras and Hopf algebroids was
considered by B\"{o}hm (\cite{Bohm04}).

Quasi-Hopf algebras have been introduced by Drinfeld (\cite{Drinfeld90}) and
have lately attracted much attention in both mathematics and physics (\cite%
{Altschuler92}, \cite{Majid95}). So it is desirable to see if it is possible
to generalize the above results also to the case of quasi-Hopf algebras.
This is the purpose of this paper.

If $H$ is a Hopf algebra, to define a Galois extension, one usually starts
with an $H$-comodule algebra $\mathcal{A}$ and its subalgebra of
coinvariants $\mathcal{A}^{coH}$. But the ordinary definition of
coinvariants does not work any more in the quasi-Hopf setting. A possible
approach was suggested in \cite{Panaite05}, but only in the case of a
morphism $H\longrightarrow \mathcal{A}$ of right comodule algebras. If we
turn to the finite dimensional case and work with module algebras,
everything seems to be fitting. So we start with $H$ a finite dimensional
quasi-Hopf algebra, $A$ a left $H$-module algebra and $B=A^{H}$ the
subalgebra of invariants, which in this case is associative, while $A$ is
not necessarily so. All the ingredients concerning quasi-Hopf actions were
already defined in \cite{Bulacu00} and \cite{BulacuPanaite00}. In order to
get a Morita context between the smash product $A\#H$ and$\ B$, we need a
right action of $A\#H$ on the linking bimodule $A$ which is not obvious. The
main ingredient in finding this action will be the \textbf{Remark} \ref{str
de modul stang}, which uses a formula for $S(t)$ (\ref{S(t)}), where $t\in H$
is a fixed nonzero left integral and $S$ is the antipode of $H$.

Next, we define Galois extensions and, as in the Hopf case, we construct two
canonical maps, showing that they are equivalent. A Frobenius-type
isomorphism (\textbf{Proposition} \ref{morfism frobenius}) is used to show
the equivalence between the canonical Galois map and one of the Morita maps.

In order to produce examples of Galois extensions, we remark first that this
definition of Galois extensions is invariant to gauge transformations. Next,
for $\mathcal{A}$ a right $H$-comodule algebra (as it was defined in \cite%
{Hausser99}), the quasi-smash product $\mathcal{A}\overline{\#}H^{\ast }$ (%
\cite{Bulacu02}) is a Galois extension of $\mathcal{A}$, as in the Hopf case.

In the last part we study the surjectivity of the second Morita map, in
connection with the notion of a total integral and the injectivity of
relative modules.

The paper is ended by an analogue of Schneider%
%TCIMACRO{\U{b4}}%
%BeginExpansion
\'{}%
%EndExpansion
s theorem in \cite{Schneider90a}.

Our proofs will follow the original proofs of \cite{Cohen90} and \cite%
{Beattie97a}. The main obstacle for the generalization is the
comultiplication of $H$ and the multiplication of $A$, which are no longer
coassociative, respectively associative. These difficulties can be overcome
by considering suitable elements that have been defined by Hausser and Nill (%
\cite{Hausser99}, \cite{Hausser99a}, \cite{HausserNill99}) and their
properties, which allow us simplify the computations.

\section{Preliminaries}

In this section we recall some definitions and results and fix notations.
Throughout the paper we work over some base field $k$. Tensor products,
algebras, linear spaces, etc. will be over $k$. Unadorned $\otimes $ means $%
\otimes _{k}$. An introduction to the study of quasi-bialgebras and
quasi-Hopf algebras can be found in \cite{Drinfeld90} or \cite{Kassel95}.

\subsection{Quasi-Hopf Algebras and their integrals}

\begin{definition}
A \textbf{quasi-bialgebra}\ means $(H,\Delta ,\varepsilon ,\phi )$ where $H$
is an associative algebra with unit, $\phi $ is an invertible element in $%
H\otimes H\otimes H$ (the associator), $\Delta :H\longrightarrow H\otimes H$
(the coproduct) and $\varepsilon :H\longrightarrow k$ (the counit) are
algebra homomorphisms, such that:%
%TCIMACRO{\TeXButton{subecuat}{\begin{subequations}}}%
%BeginExpansion
\begin{subequations}%
%EndExpansion
\begin{eqnarray}
\phi (\Delta \otimes I)\Delta (h)\phi ^{-1} &=&(I\otimes \Delta )\Delta (h)
\label{coasoc} \\
(I\otimes \varepsilon )\Delta (h) &=&(\varepsilon \otimes I)\Delta (h)=h
\label{counit} \\
(I\otimes I\otimes \Delta )(\phi )(\Delta \otimes I\otimes I)(\phi )
&=&(1\otimes \phi )(I\otimes \Delta \otimes I)(\phi )(\phi \otimes 1)
\label{pentagon} \\
(I\otimes \varepsilon \otimes I)(\phi ) &=&1\otimes 1  \label{ficounit}
\end{eqnarray}%
%TCIMACRO{\TeXButton{sfarsit}{\end{subequations}}}%
%BeginExpansion
\end{subequations}%
%EndExpansion
hold for all $h\in H$.
\end{definition}

The identities (\ref{coasoc})-(\ref{ficounit}) also imply $(\varepsilon
\otimes I\otimes I)(\phi )=(I\otimes I\otimes \varepsilon )(\phi )=1\otimes
1 $.

As for Hopf algebras, we use the Sweedler's notation $\Delta
(h)=h_{1}\otimes h_{2}$, but since $\Delta $ is only quasi-coassociative we
adopt the further convention: 
\begin{equation*}
(\Delta \otimes I)\Delta (h)=h_{1_{1}}\otimes h_{1_{2}}\otimes h_{2}\text{
and }(I\otimes \Delta )\Delta (h)=h_{1}\otimes h_{2_{1}}\otimes h_{2_{2}}
\end{equation*}%
for all $h\in H$.

We shall denote the tensor components of $\phi $ by capital letters, and
those of $\phi ^{-1}$ by small letters, namely 
\begin{eqnarray*}
\phi &=&X^{1}\otimes X^{2}\otimes X^{3}=Y^{1}\otimes Y^{2}\otimes Y^{3}=...
\\
\phi ^{-1} &=&x^{1}\otimes x^{2}\otimes x^{3}=y^{1}\otimes y^{2}\otimes
y^{3}=...
\end{eqnarray*}%
suppressing the summation symbol $\Sigma $.

\begin{definition}
A quasi-bialgebra is called a \textbf{quasi-Hopf algebra} if there is an
anti-algebra homomorphism $S:H\longrightarrow H$ (the antipode) and elements 
$\alpha $, $\beta \in H$ such that 
%TCIMACRO{\TeXButton{subecuat}{\begin{subequations}}}%
%BeginExpansion
\begin{subequations}%
%EndExpansion
\begin{eqnarray}
S(h_{1})\alpha h_{2} &=&\varepsilon (h)\alpha  \label{alfa} \\
h_{1}\beta S(h_{2}) &=&\varepsilon (h)\beta  \label{beta} \\
X^{1}\beta S(X^{2})\alpha X^{3} &=&1  \label{fibetaalfa} \\
S(x^{1})\alpha x^{2}\beta S(x^{3}) &=&1  \label{fi-1alfabeta}
\end{eqnarray}%
%TCIMACRO{\TeXButton{sfarsit}{\end{subequations}}}%
%BeginExpansion
\end{subequations}%
%EndExpansion
hold for each $h\in H$.
\end{definition}

The axioms for a quasi-Hopf algebra imply that $\varepsilon \circ
S=\varepsilon $ and $\varepsilon (\alpha )\varepsilon (\beta )=1$ so, by
rescaling $\alpha $ and $\beta $, we may assume without loss of generality
that $\varepsilon (\alpha )=\varepsilon (\beta )=1$.

In this article we consider only finite dimensional quasi-Hopf algebra $H$.
In this case the antipode of $H$ is always bijective by \cite{Bulacu03}.

Together with a quasi-Hopf algebra $H=(H,\Delta ,\varepsilon ,\phi ,S,\alpha
,\beta )$, we also have $H^{op}$, $H^{cop}$, and $H^{op,cop}$ as quasi-Hopf
algebras, where "op" means opposite multiplication and "cop" means opposite
comultiplication. The quasi-Hopf structures are obtained by putting $\phi
_{op}=\phi ^{-1}$, $\phi _{cop}=(\phi ^{-1})^{321}$, $\phi _{op,cop}=\phi
^{321}$, $S_{op}=S_{cop}=(S_{op,cop})^{-1}=S^{-1}$, $\alpha
_{op}=S^{-1}(\beta )$,$\,\ \alpha _{cop}=S^{-1}(\alpha )$, $\alpha
_{op,cop}=\beta $, $\beta _{op}=S^{-1}(\alpha )$, $\beta _{cop}=S^{-1}(\beta
)$ and $\beta _{op,cop}=\alpha $.

Suppose that $(H,\Delta ,\varepsilon ,\phi )$ is a quasi-bialgebra. Then the
category $_{H}\mathcal{M}$ of left $H$-modules is monoidal in the following
way: for $U,V\in $ $_{H}\mathcal{M}$ the tensor product $U\otimes V$ is an $%
H $-module by $h(u\otimes v)=h_{1}u\otimes h_{2}v$. The base field $k$ is an 
$H $-module via $\varepsilon $. The canonical morphisms $U\simeq U\otimes
k\simeq k\otimes U$ are $H$-linear for $U\in $ $_{H}\mathcal{M}$. The map 
\begin{eqnarray*}
\Phi _{U,V,W} &:&(U\otimes V)\otimes W\longrightarrow U\otimes (V\otimes W)
\\
u\otimes v\otimes w &\longrightarrow &X^{1}u\otimes X^{2}v\otimes X^{3}w
\end{eqnarray*}%
for $U,V,W\in $ $_{H}\mathcal{M}$ is $H$-linear as a consequence of (\ref%
{coasoc}), and satisfies Mac Lane's pentagon axiom for a monoidal category
as a consequence of (\ref{pentagon}).

In the Hopf algebra case, the antipode is an anti-coalgebra map. In order to
have a similar property in the quasi-Hopf setting, Drinfeld (\cite%
{Drinfeld90}) introduced a gauge element $f\in H\otimes H$, which obeys the
following:%
\begin{eqnarray}
f\Delta (h)f^{(-1)} &=&(S\otimes S)\Delta ^{cop}S^{-1}(h)  \label{fdeltaf-1}
\\
(S\otimes S\otimes S)(\phi ^{321}) &=&(1\otimes f)(I\otimes \Delta )(f)\phi
(\Delta \otimes I)(f^{(-1)})(f^{(-1)}\otimes 1)  \label{fif} \\
(I\otimes \varepsilon )(f) &=&(\varepsilon \otimes I)(f)=1  \label{gauge}
\end{eqnarray}%
for all $h\in H$.

Following \cite{Hausser99}, \cite{Hausser99a}, \cite{HausserNill99},we may
define the elements%
\begin{eqnarray}
p_{L} &=&X^{2}S^{-1}(X^{1}\beta )\otimes X^{3}  \label{formpl} \\
q_{L} &=&S(x^{1})\alpha x^{2}\otimes x^{3}  \label{formql} \\
p_{R} &=&x^{1}\otimes x^{2}\beta S(x^{3})  \label{formpr} \\
q_{R} &=&X^{1}\otimes S^{-1}(\alpha X^{3})X^{2}  \label{formqr}
\end{eqnarray}%
One may note that in $H^{op}$ the roles of $p_{R}$ and $q_{R}$ (respectively 
$p_{L}$ and $q_{L}$) interchange, and in $H^{cop}$ one is passing from $%
p_{L} $ to $p_{R}$\ (respectively from $q_{L}$ to $q_{R}$), so whenever we
have a relation concerning these elements, it is enough to pass to $H^{op}$
or to $H^{cop}$ to get the other three relations. We shall state here only
the relations used in this paper. Following \cite{Hausser99}, they obey the
following 
\begin{eqnarray}
\Delta (h_{2})p_{L}(S^{-1}(h_{1})\otimes 1) &=&p_{L}(1\otimes h)  \label{pl}
\\
(S(h_{1})\otimes 1)q_{L}\Delta (h_{2}) &=&(1\otimes h)q_{L}  \label{ql} \\
\Delta (h_{1})p_{R}(1\otimes S(h_{2})) &=&p_{R}(h\otimes 1)  \label{pr} \\
(1\otimes S^{-1}(h_{2}))q_{R}\Delta (h_{1}) &=&(h\otimes 1)q_{R}  \label{qr}
\end{eqnarray}%
for all $h\in H$ and%
\begin{eqnarray}
\Delta (q_{L}^{2})p_{L}(S^{-1}(q_{L}^{1})\otimes 1) &=&1\otimes 1
\label{qlpl} \\
(S(p_{L}^{1})\otimes 1)q_{L}\Delta (p_{L}^{2}) &=&1\otimes 1  \label{plql} \\
\Delta (q_{R}^{1})p_{R}(1\otimes S(q_{R}^{2})) &=&1\otimes 1  \label{qrpr} \\
(1\otimes S^{-1}(p_{R}^{2}))q_{R}\Delta (p_{R}^{1}) &=&1\otimes 1
\label{prqr}
\end{eqnarray}%
\begin{eqnarray}
\phi ^{-1}(I\otimes \Delta )(p_{L}) &=&(\Delta (X^{2})p_{L}\otimes
X^{3})(S^{-1}(X^{1})\otimes 1\otimes 1)  \label{fipl} \\
(I\otimes \Delta )(q_{L})\phi &=&(S(x^{1})\otimes 1\otimes 1)(q_{L}\Delta
(x^{2})\otimes x^{3})  \label{qlfi} \\
\phi (\Delta \otimes I)(p_{R}) &=&(x^{1}\otimes \Delta
(x^{2})p_{R})(1\otimes 1\otimes S(x^{3}))  \label{fipr} \\
(\Delta \otimes I)(q_{R})\phi ^{-1} &=&(1\otimes 1\otimes
S^{-1}(X^{3}))(X^{1}\otimes q_{R}\Delta (X^{2}))  \label{qrfi}
\end{eqnarray}

We shall need also the following elements:%
\begin{eqnarray}
V_{L} &=&(S\otimes S)(p_{L}^{21})f  \label{VL} \\
U_{L} &=&(S^{-1}\otimes S^{-1})(q_{L}f^{(-1)})^{21}  \label{UL} \\
V_{R} &=&(S^{-1}\otimes S^{-1})(fp_{R})^{21}  \label{VR} \\
U_{R} &=&f^{-1}(S\otimes S)(q_{R}^{21})  \label{UR}
\end{eqnarray}%
$U_{R}$ and $V_{R}$ were introduced by Nill and Hausser in \cite%
{HausserNill99}, and $U_{L}$ and $V_{L}$ are their analogues by passing from 
$H$ to $H^{cop}$, as explained above. They obey the following relations:%
\begin{eqnarray}
(1\otimes h_{1})V_{L}\Delta S(h_{2}) &=&(S(h)\otimes 1)V_{L} \\
\Delta S^{-1}(h_{2})U_{L}(1\otimes h_{1}) &=&U_{L}(S^{-1}(h)\otimes 1) \\
(h_{2}\otimes 1)V_{R}\Delta S^{-1}(h_{1}) &=&(1\otimes S^{-1}(h))V_{R} \\
\Delta S(h_{1})U_{R}(h_{2}\otimes 1) &=&U_{R}(1\otimes S(h))
\end{eqnarray}

We need also to notice that $H$ finite dimensional implies $H^{\ast }$ is a
coassociative coalgebra, with coproduct given by%
\begin{equation*}
\Delta ^{\ast }(h^{\ast })=h^{\ast }\circ \mu _{H}\text{, }\forall h^{\ast
}\in H^{\ast }
\end{equation*}%
The comultiplication $\Delta $ on $H$ allows us to define a multiplication
on $H^{\ast }$ (the convolution product) by 
\begin{equation*}
(h^{\ast }g^{\ast })(h)=h^{\ast }(h_{1})g^{\ast }(h_{2})\text{, }\forall
h^{\ast }\in H^{\ast }\text{, }g,h\in H
\end{equation*}%
But this is no longer associative. In fact $H^{\ast }$, endowed with this
multiplication, is an algebra in the monoidal category of $H$-bimodules and
we have that 
\begin{equation*}
(h^{\ast }g^{\ast })l^{\ast }=(X^{1}\rightharpoonup h^{\ast }\leftharpoonup
x^{1})[(X^{2}\leftharpoonup g^{\ast }\leftharpoonup
x^{2})(X^{3}\rightharpoonup l^{\ast }\rightharpoonup x^{3})]
\end{equation*}%
where $h^{\ast },g^{\ast },l^{\ast }\in H^{\ast }$. By $\leftharpoonup $ and 
$\rightharpoonup $ we denote the usual right, respectively left $H$-action
on $H^{\ast }$:%
\begin{equation*}
(h\rightharpoonup h^{\ast })(g)=h^{\ast }(gh)\text{, }(h^{\ast
}\leftharpoonup h)(g)=h^{\ast }(hg)
\end{equation*}%
for all $h^{\ast }\in H^{\ast }$, $g,h\in H$.

Although $H^{\ast }$ is not an algebra, we still keep the notation of the
Hopf case for the weak action of $H^{\ast }$ on $H$:%
\begin{equation*}
h^{\ast }\rightharpoonup h=h^{\ast }(h_{2})h_{1}\text{, }h\leftharpoonup
h^{\ast }=h^{\ast }(h_{1})h_{2}
\end{equation*}%
for all $h^{\ast }\in H^{\ast }$, $h\in H$.

We denote by $\int_{H}^{l}$ the space of left integrals in $H$ and by $t\in
H $ a nonzero left integral. As $\int_{H}^{l}$ is an ideal of $H$ (which is
one dimensional by \cite{HausserNill99}), there is only one algebra morphism
(the modular element) $\gamma \in H^{\ast }$ which satisfies 
\begin{equation}
th=\gamma (h)t  \label{gama}
\end{equation}%
for all $h\in H$. Denote by $\Lambda =\gamma (q_{L}^{2})q_{L}^{1}$. The next
result was proved by Bulacu and Caenepeel in \cite{Bulacu03}, \textbf{%
Proposition} 4.9 for a dual quasi-Hopf algebra. But if $H$ is finite
dimensional, $H^{\ast }$ is a dual quasi-Hopf algebra, so we can restate
their result in terms of quasi-Hopf algebras to get the following:

\begin{proposition}
\label{S(t)}Let $H$ be a finite dimensional quasi-Hopf algebra. With the
above notation we have 
\begin{equation*}
S(t)=\Lambda (\gamma \rightharpoonup t).
\end{equation*}
\end{proposition}

The following statement\textbf{\ }can be obtained from the same authors
quoted above (\cite{Bulacu03}, \textbf{Theorem} 2.2) by passing from $H$ to $%
H^{cop}$:

\begin{theorem}
Let $H$ be a finite dimensional quasi-Hopf algebra, $(e_{i})_{i=1,n}$ a
basis of $H$ and $(e^{i})_{i=1,n}$ the dual basis of $H^{\ast }$. We define
the map%
\begin{equation*}
P:H\longrightarrow H\text{, }P(h)=%
\sum_{i=1}^{n}e^{i}(S^{-2}(q_{L}^{1}e_{i1}S(\beta ))h)q_{L}^{2}e_{i2}
\end{equation*}%
for all $h\in H$. Then:

\begin{enumerate}
\item $P(h)\in \int_{H}^{l}$, for all $h\in H$ and there is $h\in H$ such
that $P(h)\neq 0$;

\item The map $\Theta :\int_{H}^{l}\otimes H^{\ast }\longrightarrow H$, 
\begin{equation*}
\Theta (t\otimes h^{\ast })=h^{\ast
}(q_{L}^{1}t_{1}p_{L}^{1})q_{L}^{2}t_{2}p_{L}^{2}\text{, }\forall t\in
\int_{H}^{l}\text{, }h^{\ast }\in H^{\ast }
\end{equation*}%
is an isomorphism of left $H$-modules with inverse given by 
\begin{equation*}
\Theta ^{-1}(h)=\sum_{i=1}^{n}P(e_{i}h)\otimes e^{i}S^{-1}\text{, }\forall
h\in H,
\end{equation*}%
where $\int_{H}^{l}\otimes H^{\ast }$ is a left $H$-module via $h(t\otimes
h^{\ast })=t\otimes h^{\ast }\leftharpoonup S(h)$, $\forall h\in H$, $t\in
\int_{H}^{l}$, $h^{\ast }\in H^{\ast }$.
\end{enumerate}
\end{theorem}

This allows us to consider the following Frobenius-type isomorphism:

\begin{corollary}
\label{morfism frobenius}Let $H$ be a finite dimensional quasi-Hopf algebra
and $t\in \int_{H}^{l}$ a nonzero left integral. Then the map 
\begin{equation*}
\theta _{t}:H^{\ast }\longrightarrow H\text{, }\theta _{t}(h^{\ast
})=h^{\ast }(q_{L}^{1}t_{1}p_{L}^{1})q_{L}^{2}t_{2}p_{L}^{2}\text{, }\forall
h^{\ast }\in H^{\ast }
\end{equation*}%
is a left $H$-module isomorphism.
\end{corollary}

By \cite{Bulacu03} or \cite{Mason05}, we have, for $t\in H$ a left integral
and all $h\in H$:

\begin{eqnarray}
(S(h)\otimes 1)q_{L}\Delta (t) &=&(1\otimes h)q_{L}\Delta (t)  \label{Sqlt}
\\
(S(h)\otimes 1)q_{R}\Delta (t) &=&(1\otimes h)q_{R}\Delta (t)  \label{Sqrt}
\end{eqnarray}%
and

\begin{eqnarray}
\Delta (t) &=&(\beta \otimes 1)q_{L}\Delta (t)=(\beta \otimes 1)q_{R}\Delta
(t)  \label{betaql} \\
&=&(1\otimes S^{-1}(\beta ))q_{L}\Delta (t)=(1\otimes S^{-1}(\beta
))q_{R}\Delta (t)  \label{sbetaql}
\end{eqnarray}

\subsection{$H$-module algebras}

Recall from \cite{BulacuPanaite00} the notion of a module algebra over a
quasi-bialgebra.

\begin{definition}
Let $H$ be a quasi-bialgebra and $A$ a linear space. We say that $A$ is a
(left) $H$-\textbf{module algebra} if $A$ is an algebra in the monoidal
category $_{H}\mathcal{M}$, i.e. $A$ is a left $H$-module which has a
multiplication and a usual unit $1_{A}$ satisfying the following conditions: 
%TCIMACRO{\TeXButton{subecuat}{\begin{subequations}}}%
%BeginExpansion
\begin{subequations}%
%EndExpansion
\begin{eqnarray}
(ab)c &=&(X^{1}\cdot a)[(X^{2}\cdot b)(X^{3}\cdot c)]  \label{smashasoc} \\
h\cdot (ab) &=&(h_{1}\cdot a)(h_{2}\cdot b)  \label{smashmod} \\
h\cdot 1_{A} &=&\varepsilon (h)1_{A}  \label{2.11}
\end{eqnarray}%
%TCIMACRO{\TeXButton{sfarsit}{\end{subequations}}}%
%BeginExpansion
\end{subequations}%
%EndExpansion
for all $a,b,c\in A$ and $h\in H$, where $h\otimes a\longrightarrow h\cdot a$
is the $H$-module structure \ of $A$.
\end{definition}

For an $H$-module algebra $A$, we may define the \textbf{smash product} $%
A\#H $ as in \cite{BulacuPanaite00}: as a vector space $A\#H$ is $A\otimes H$
with multiplication given by 
\begin{equation}
(a\#h)(b\#g)=(x^{1}\cdot a)(x^{2}h_{1}\cdot b)\#x^{3}h_{2}g
\label{smashproduct}
\end{equation}%
for all $a,b\in A$ and $g,h\in H$. Then $A\#H$ becomes an associative
algebra with unit $1_{A}\#1$. Also for the $H$-module algebra $A$ we may
define the subalgebra of invariants $B=A^{H}$, that is%
\begin{equation*}
B=\left\{ a\in A\mid h\cdot a=\varepsilon (h)a,\forall h\in H\right\}
\end{equation*}%
Remark that $B$ is an associative algebra, and $A$ is a left and right $B$%
-module in a natural way. Also, on $A$ we have a left $A\#H$-module
structure given by 
\begin{equation*}
(a\#h)b=a(h\cdot b)
\end{equation*}%
for all $a,b\in A$ and $h\in H$.

As $A$ is an algebra in the monoidal category of left $H$-modules, it is
natural to consider modules over $A$ (left or right) in this category. These
were called relative Hopf modules in \cite{BulacuPanaite00}:

\begin{definition}
Let $H$ be a quasi-bialgebra and $A$ a left $H$-module algebra. A $k$-vector
space $M$ is called a left $(H,A)$-\textbf{Hopf module} if $M$ is a left $H$%
-module and also a left $A$-module in the monoidal category $_{H}\mathcal{M}$%
, i.e. $A$ acts on $M$ to the left with action denoted $a\otimes
m\longrightarrow am$ such that 
%TCIMACRO{\TeXButton{subecuat}{\begin{subequations}}}%
%BeginExpansion
\begin{subequations}%
%EndExpansion
\begin{eqnarray}
(ab)m &=&(X^{1}\cdot a)[(X^{2}\cdot b)(X^{3}m)] \\
h(am) &=&(h_{1}\cdot a)(h_{2}m) \\
1_{A}m &=&m
\end{eqnarray}%
%TCIMACRO{\TeXButton{sfarsit}{\end{subequations}}}%
%BeginExpansion
\end{subequations}%
%EndExpansion
for all $h\in H$, $a,b\in A$ and $m\in M$.
\end{definition}

The category of left $(H,A)$-Hopf modules with morphisms that are left $H$%
-linear and preserve the weak $A$-action will be denoted by $_{A}(_{H}%
\mathcal{M)}$.

In \cite{BulacuPanaite00} it is proved that the categories $_{A}(_{H}%
\mathcal{M)}$ and $_{A\#H}\mathcal{M}$ are isomorphic: if $M\in $ $_{A}(_{H}%
\mathcal{M)}$, then $M\in $ $_{A\#H}\mathcal{M}$ by $(a\#h)m=a(hm)$, and if $%
M\in $ $_{A\#H}\mathcal{M}$ then $M\in $ $_{A}(_{H}\mathcal{M)}$ by $%
am=(a\#1)m$ and $hm=(1_{A}\#h)m$, where $h\in H$, $a\in A$ and $m\in M$.

Now for each $M\in $ $_{A}(_{H}\mathcal{M)}$, denote $M^{H}=\{m\in M\mathcal{%
\mid }hm=\mathcal{\varepsilon (}h)m,\forall h\in H\}$. Then $M^{H}$ becomes
naturally a left $B$-module, so we get a functor 
\begin{equation*}
_{A}(_{H}\mathcal{M)}\overset{(-)^{H}}{\mathcal{\longrightarrow }}\text{ }%
_{B}\mathcal{M}
\end{equation*}%
Conversely, for $M\in $ $_{B}\mathcal{M}$, we have $A\otimes _{B}M\in $ $%
_{A}(_{H}\mathcal{M)}$ by 
\begin{eqnarray*}
h(a\otimes _{B}m) &=&h\cdot a\otimes _{B}m \\
b(a\otimes _{B}m) &=&ba\otimes _{B}m
\end{eqnarray*}%
As in the classical Hopf algebra case, we obtain the following:

\begin{proposition}
\label{functorul indus}The induced functor $A\otimes _{B}(-)$ is a left
adjoint for the functor of invariants $(-)^{H}$%
\begin{equation*}
_{B}\mathcal{M}\overset{A\otimes _{B}(-)}{\underset{(-)^{H}}{\mathcal{%
\rightleftarrows }}}\text{ }_{A}(_{H}\mathcal{M)}
\end{equation*}
\end{proposition}

\begin{proof}
The adjunction morphisms are:%
\begin{eqnarray*}
f\in Hom_{_{A}(_{H}\mathcal{M)}}(A\mathcal{\otimes }_{B}M,N) &&\underset{%
\beta }{\overset{\alpha }{\rightleftarrows }}Hom_{B}(M,N^{H})\ni g \\
\alpha (f)(m) &=&f(1_{A}\otimes _{B}m) \\
\beta (g)(a\otimes _{B}m) &=&ag(m)
\end{eqnarray*}%
for all $m\in M$, $a\in A$.
\end{proof}

\begin{remark}
As $_{A}(_{H}\mathcal{M)}$ and $_{A\#H}\mathcal{M}$ are isomorphic, we get
also the adjunction $_{B}\mathcal{M}\overset{A\otimes _{B}(-)}{\underset{%
(-)^{H}}{\mathcal{\rightleftarrows }}}$ $_{A\#H}\mathcal{M}$.
\end{remark}

\section{The Morita context}

\subsection{Construction}

In this section we construct a Morita context connecting $A\#H$ and $B=A^{H}$%
, where $A$ is our left $H$-module algebra. By \cite{BulacuPanaite00}, 
\textbf{Proposition} 2.7, the map $H\longrightarrow A\#H$, $h\longrightarrow
1_{A}\#h$ is an algebra morphism. It induces a structure of left $H$-module
on $A\#H$ by 
\begin{equation}
h\times (a\#g)=(1_{A}\#h)(a\#g)=h_{1}\cdot a\#h_{2}g  \label{hmodsmash}
\end{equation}%
for all $a\in A$ and $g,h\in H$. Also, the tensor product $H\otimes A$ has
the structure of a left $H$-module with action induced by the multiplication
of $H$:%
\begin{equation*}
h\ast (g\otimes a)=hg\otimes a
\end{equation*}%
for all $a\in A$ and $g,h\in H$.

\begin{proposition}
\label{isomHtensorAcuAtensorH}Under the above assumptions, the map $\varphi
:H\otimes A\longrightarrow A\#H$, $\varphi (h\otimes a)=h_{1}p_{L}^{1}\cdot
a\#h_{2}p_{L}^{2}$ is an isomorphism of left $H$-modules, with inverse $%
\varphi ^{-1}(a\#h)=q_{L}^{2}h_{2}\otimes S^{-1}(q_{L}^{1}h_{1})\cdot a$.
\end{proposition}

\begin{proof}
We have 
\begin{eqnarray*}
\varphi (h\ast (g\otimes a)) &=&\varphi (hg\otimes a) \\
&=&h_{1}g_{1}p_{L}^{1}\cdot a\#h_{2}g_{2}p_{L}^{2} \\
(\ref{hmodsmash}) &=&h\times (g_{1}p_{L}^{1}\cdot a\#g_{2}p_{L}^{2}) \\
&=&h\times \varphi (g\otimes a)
\end{eqnarray*}%
for all $a\in A$ and $g,h\in H$, so $\varphi $ is $H$-linear. Let's check
that $\varphi $ and $\varphi ^{-1}$ are inverses to each other:%
\begin{eqnarray*}
\varphi \varphi ^{-1}(a\#h) &=&\varphi (q_{L}^{2}h_{2}\otimes
S^{-1}(q_{L}^{1}h_{1})\cdot a) \\
&=&q_{L1}^{2}h_{2_{1}}p_{L}^{1}S^{-1}(q_{L}^{1}h_{1})\cdot
a\#q_{L2}^{2}h_{2_{2}}p_{L}^{2} \\
(\ref{pl}) &=&q_{L1}^{2}p_{L}^{1}S^{-1}(q_{L}^{1})\cdot
a\#q_{L2}^{2}p_{L}^{2}h \\
(\ref{qlpl}) &=&a\#h
\end{eqnarray*}%
and 
\begin{eqnarray*}
\varphi ^{-1}\varphi (h\otimes a) &=&\varphi ^{-1}(h_{1}p_{L}^{1}\cdot
a\#h_{2}p_{L}^{2}) \\
&=&q_{L}^{2}h_{2_{2}}p_{L2}^{2}\otimes
S^{-1}(q_{L}^{1}h_{2_{1}}p_{L1}^{2})h_{1}p_{L}^{1}\cdot a \\
(\ref{ql}) &=&hq_{L}^{2}p_{L2}^{2}\otimes
S^{-1}(q_{L}^{1}p_{L1}^{2})p_{L}^{1}\cdot a \\
(\ref{plql}) &=&h\otimes a
\end{eqnarray*}%
for all $a\in A$ and $h\in H$.
\end{proof}

\begin{remark}
\label{str de modul stang}Restricting the above isomorphism, we get an
isomorphism of $H$-invariants: $(A\#H)^{H}\simeq (H\otimes A)^{H}$. But $H$
is acting on $H\otimes A$ by left multiplication on the first component, so $%
(H\otimes A)^{H}=H^{H}\otimes A=\int_{H}^{l}\otimes A\simeq A$.\ On the
other side, $A\#H$ is a right $A\#H$-module with action induced by
multiplication, and the restriction $(A\#H)^{H}$ remains a right $A\#H$%
-module, because for $a\#h\in (A\#H)^{H}$, $b\#g\in A\#H$ and $l\in H$ we
have%
\begin{eqnarray*}
l\times ((a\#h)(b\#g)) &=&(1\#l)[(a\#h)(b\#g)] \\
&=&[(1\#l)(a\#h)](b\#g) \\
&=&\varepsilon (l)(a\#h)(b\#g)
\end{eqnarray*}%
Hence, the isomorphism of \textbf{Proposition} \ref{isomHtensorAcuAtensorH}
induces a structure of right $A\#H$-module on $A\simeq \int_{H}^{l}\otimes A$%
. Explicitly, this means%
%TCIMACRO{\TeXButton{breakdisplay}{\begin{allowdisplaybreaks}}}%
%BeginExpansion
\begin{allowdisplaybreaks}%
%EndExpansion
\begin{eqnarray*}
a\otimes (b\#h) &\longrightarrow &(t\otimes a)\otimes (b\#h)\overset{\varphi
\otimes I}{\longrightarrow }(t_{1}p_{L}^{1}\cdot a\#t_{2}p_{L}^{2})\otimes
(b\#h) \\
&\longrightarrow &(t_{1}p_{L}^{1}\cdot a\#t_{2}p_{L}^{2})(b\#h) \\
&=&(x^{1}t_{1}p_{L}^{1}\cdot a)(x^{2}t_{2_{1}}p_{L1}^{2}\cdot
b)\#x^{3}t_{2_{2}}p_{L2}^{2}h \\
(\ref{coasoc}) &=&(t_{1_{1}}x^{1}p_{L}^{1}\cdot
a)(t_{1_{2}}x^{2}p_{L1}^{2}\cdot b)\#t_{2}x^{3}p_{L2}^{2}h \\
(\ref{fipl}) &=&t_{1}X^{2}\cdot \lbrack (p_{L}^{1}S^{-1}(X^{1})\cdot
a)(p_{L}^{2}\cdot b)]\#t_{2}X^{3}h \\
&&\overset{\varphi ^{-1}}{\longrightarrow }q_{L}^{2}t_{2_{2}}X_{2}^{3}h_{2}%
\otimes S^{-1}(q_{L}^{1}t_{2_{1}}X_{1}^{3}h_{1})t_{1}X^{2} \\
&&\cdot \lbrack (p_{L}^{1}S^{-1}(X^{1})\cdot a)(p_{L}^{2}\cdot b)] \\
(\ref{ql}) &=&tq_{L}^{2}X_{2}^{3}h_{2}\otimes
S^{-1}(q_{L}^{1}X_{1}^{3}h_{1})X^{2}\cdot \lbrack
(p_{L}^{1}S^{-1}(X^{1})\cdot a)(p_{L}^{2}\cdot b)] \\
(\ref{gama}) &=&t\otimes S^{-1}(q_{L}^{1}X_{1}^{3}h_{1})X^{2}\cdot \lbrack
(p_{L}^{1}S^{-1}(X^{1})\cdot a)(p_{L}^{2}\cdot b)]\gamma
(q_{L}^{2}X_{2}^{3}h_{2}) \\
&=&t\otimes S^{-1}(S(X^{2})\Lambda (\gamma \rightharpoonup X^{3}h))\cdot 
\left[ (p_{L}^{1}S^{-1}(X^{1})\cdot a)(p_{L}^{2}\cdot b)\right]
\end{eqnarray*}%
%TCIMACRO{\TeXButton{endbreakdisplay}{\end{allowdisplaybreaks}}}%
%BeginExpansion
\end{allowdisplaybreaks}%
%EndExpansion
where we denoted $\Lambda =q_{L}^{1}\gamma (q_{L}^{2})$. Therefore, the
right $A\#H$-module structure on $A$ is given by%
\begin{equation}
a\cdot _{\gamma }(b\#h)=S^{-1}(S(X^{2})\Lambda (\gamma \rightharpoonup
X^{3}h))\cdot \left[ (p_{L}^{1}S^{-1}(X^{1})\cdot a)(p_{L}^{2}\cdot b)\right]
\label{str la dr de a smash h mod}
\end{equation}%
This can be also checked directly, but very long calculations are involved.
\end{remark}

By (\ref{str la dr de a smash h mod}), $A$ becomes a $(B,A\#H)$-bimodule. In
the previous section we endowed $A$ with an $(A\#H,B)$-bimodule structure,
so now we may pass to the next step:

\begin{theorem}
Consider the maps: 
\begin{eqnarray}
\left( -,-\right) &:&A\otimes _{A\#H}A\longrightarrow B\text{, }(a,b)=t\cdot %
\left[ (p_{L}^{1}\cdot a)(p_{L}^{2}\cdot b)\right]  \label{moritaB} \\
\left[ -,-\right] &:&A\otimes _{B}A\longrightarrow A\#H\text{, }\left[ a,b%
\right] =(a\#t)(p_{L}^{1}\cdot b\#p_{L}^{2})  \label{moritasmash}
\end{eqnarray}%
Then $(A\#H$, $B$, $_{A\#H}A_{B}$, $_{B}A_{A\#H}$, $\left( -,-\right) $, $%
\left[ -,-\right] )$ is a Morita context.
\end{theorem}

\begin{proof}
We check only the conditions of a Morita context which are more difficult
because of the $A\#H$-module structures involved, the others are left to the
reader:

\begin{itemize}
\item $\left( -,-\right) $ is well-defined ($A\#H$-balanced):%
\begin{eqnarray*}
(a\cdot _{\gamma }(b\#h),c) &=&(S^{-1}(S(X^{2})\Lambda (\gamma
\rightharpoonup X^{3}h))\cdot \lbrack (p_{L}^{1}S^{-1}(X^{1})\cdot
a)(p_{L}^{2}\cdot b)],c) \\
(\ref{moritaB}) &=&t\cdot \{(P_{L}^{1}S^{-1}(S(X^{2})\Lambda (\gamma
\rightharpoonup X^{3}h))\cdot \left[ (p_{L}^{1}S^{-1}(X^{1})\cdot
a)(p_{L}^{2}\cdot b)\right] )(P_{L}^{2}\cdot c)\} \\
(\ref{gama}) &=&tq_{L}^{2}X_{2}^{3}h_{2}\cdot
\{(P_{L}^{1}S^{-1}(q_{L}^{1}X_{1}^{3}h_{1})X^{2}\cdot \left[
(p_{L}^{1}S^{-1}(X^{1})\cdot a)(p_{L}^{2}\cdot b)\right] )(P_{L}^{2}\cdot
c)\} \\
(\ref{pl}) &=&tq_{L}^{2}\cdot \{(P_{L}^{1}S^{-1}(q_{L}^{1})X^{2}\cdot \left[
(p_{L}^{1}S^{-1}(X^{1})\cdot a)(p_{L}^{2}\cdot b)\right] )(P_{L}^{2}X^{3}h%
\cdot c)\} \\
(\ref{qlpl}) &=&t\cdot \left\{ (X^{2}\cdot \left[ (p_{L}^{1}S^{-1}(X^{1})%
\cdot a)(p_{L}^{2}\cdot b)\right] )(X^{3}h\cdot c)\right\} \\
(\ref{fipl}) &=&t\cdot \left\{ \left[ (x^{1}p_{L}^{1}\cdot
a)(x^{2}p_{L1}^{2}\cdot b)\right] (x^{3}p_{L2}^{2}h\cdot c)\right\} \\
(\ref{smashasoc}),(\ref{smashmod}) &=&t\cdot \left\{ (p_{L}^{1}\cdot
a)[p_{L}^{2}\cdot (b(h\cdot c))\right\} \\
(\ref{moritaB}) &=&(a,b(h\cdot c)) \\
&=&(a,(b\#h)c)
\end{eqnarray*}%
for each $a,b,c\in A$ and $h\in H$, where $P_{L}^{1}\otimes P_{L}^{2}$ is
another copy of $p_{L}$.

\item $\left[ -,-\right] $ is $A\#H$-bilinear: for the left linearity, we
compute that%
\begin{eqnarray*}
\left[ (a\#h)b,c\right] &=&\left[ a(hb),c\right] \\
(\ref{moritasmash}) &=&(a(hb)\#t)(p_{L}^{1}c\#p_{L}^{2}) \\
&=&\left[ (x^{1}a)(x^{2}h_{1}b)\#x^{3}h_{2}t\right] (p_{L}^{1}c\#p_{L}^{2})
\\
(\ref{smashproduct}) &=&(a\#h)(b\#t)(p_{L}^{1}c\#p_{L}^{2}) \\
&=&(a\#h)\left[ b,c\right]
\end{eqnarray*}%
where in the third line we used the fact that $t$ is a left integral. For
the right $A\#H$-linearity it's a little more complicated because of the
right action of $A\#H$ on $A$:%
\begin{eqnarray*}
\left[ a,b\cdot _{\gamma }(c\#h)\right] &=&\left[ a,S^{-1}(S(X^{2})\Lambda
(\gamma \rightharpoonup X^{3}h))\cdot \left[ (p_{L}^{1}S^{-1}(X^{1})\cdot
b)(p_{L}^{2}\cdot c)\right] \right] \\
(\ref{moritasmash}) &=&(a\#t)(P_{L}^{1}S^{-1}(S(X^{2})\Lambda (\gamma
\rightharpoonup X^{3}h))\cdot \left[ (p_{L}^{1}S^{-1}(X^{1})\cdot
b)(p_{L}^{2}\cdot c)\right] \#P_{L}^{2}) \\
(\ref{gama})
&=&(a%
\#tq_{L}^{2}X_{2}^{3}h_{2})(P_{L}^{1}S^{-1}(q_{L}^{1}X_{1}^{3}h_{1})X^{2}%
\cdot \left[ (p_{L}^{1}S^{-1}(X^{1})\cdot b)(p_{L}^{2}\cdot c)\right]
\#P_{L}^{2}) \\
(\ref{smashproduct}) &=&(x^{1}\cdot
a)(x^{2}t_{1}q_{L1}^{2}X_{2_{1}}^{3}h_{2_{1}}P_{L}^{1}S^{-1}(q_{L}^{1}X_{1}^{3}h_{1})X^{2}\cdot %
\left[ (p_{L}^{1}S^{-1}(X^{1})\cdot b)(p_{L}^{2}\cdot c)\right] ) \\
&&\#x^{3}t_{2}q_{L2}^{2}X_{2_{2}}^{3}h_{2_{1}}P_{L}^{2} \\
(\ref{pl}) &=&(x^{1}\cdot
a)(x^{2}t_{1}q_{L1}^{2}P_{L}^{1}S^{-1}(q_{L}^{1})X^{2}\cdot \left[
(p_{L}^{1}S^{-1}(X^{1})\cdot b)(p_{L}^{2}\cdot c)\right] ) \\
&&\#x^{3}t_{2}q_{L2}^{2}P_{L}^{2}X^{3}h \\
(\ref{qlpl}) &=&(x^{1}\cdot a)(x^{2}t_{1}X^{2}\cdot \left[
(p_{L}^{1}S^{-1}(X^{1})\cdot b)(p_{L}^{2}\cdot c)\right] )\#x^{3}t_{2}X^{3}h
\\
(\ref{fipl}) &=&(x^{1}\cdot a)(x^{2}t_{1}\cdot \left[ (y^{1}p_{L}^{1}\cdot
b)(y^{2}p_{L1}^{2}\cdot c)\right] \#x^{3}t_{2}y^{3}p_{L2}^{2}h \\
(\ref{smashproduct}) &=&(a\#t)((y^{1}p_{L}^{1}\cdot b)(y^{2}p_{L1}^{2}\cdot
c)\#y^{3}p_{L2}^{2}h) \\
(\ref{smashproduct}) &=&(a\#t)(p_{L}^{1}\cdot b\#p_{L}^{2})(c\#h) \\
(\ref{moritasmash}) &=&\left[ a,b\right] (c\#h)
\end{eqnarray*}%
for all $a,b,c\in A$ and $h\in H$, where $P_{L}^{1}\otimes P_{L}^{2}$ is
another copy of $p_{L}$.

\item associativity of the Morita map:%
\begin{eqnarray*}
a\cdot _{\gamma }\left[ b,c\right] &=&a\cdot _{\gamma }\left[
(b\#t)(p_{L}^{1}\cdot c\#p_{L}^{2})\right] \\
&=&\left[ a\cdot _{\gamma }(b\#t)\right] \cdot _{\gamma }(p_{L}^{1}\cdot
c\#p_{L}^{2}) \\
(\ref{str la dr de a smash h mod}) &=&\left\{ S^{-1}(S(X^{2})\Lambda (\gamma
\rightharpoonup X^{3}t))\cdot \left[ (P_{L}^{1}S^{-1}(X^{1})\cdot
a)(P_{L}^{2}\cdot b)\right] \right\} \cdot _{\gamma }(p_{L}^{1}c\#p_{L}^{2})
\\
&=&\left\{ S^{-1}(\Lambda (\gamma \rightharpoonup t))\cdot \left[
(P_{L}^{1}\cdot a)(P_{L}^{2}\cdot b)\right] \right\} \cdot _{\gamma
}(p_{L}^{1}c\#p_{L}^{2}) \\
(\ref{S(t)}) &=&\left\{ S^{-1}(S(t))\cdot \left[ (P_{L}^{1}\cdot
a)(P_{L}^{2}\cdot b)\right] \right\} \cdot _{\gamma }(p_{L}^{1}\cdot
c\#p_{L}^{2}) \\
&=&\left\{ t\cdot \left[ (P_{L}^{1}\cdot a)(P_{L}^{2}\cdot b)\right]
\right\} \cdot _{\gamma }(p_{L}^{1}\cdot c\#p_{L}^{2}) \\
&=&(a,b)\cdot _{\gamma }(p_{L}^{1}\cdot c\#p_{L}^{2}) \\
(\ref{str la dr de a smash h mod}) &=&S^{-1}(S(X^{2})\Lambda (\gamma
\rightharpoonup X^{3}p_{L}^{2}))\cdot \left[ (P_{L}^{1}S^{-1}(X^{1})\cdot
(a,b))(P_{L}^{2}p_{L}^{1}\cdot c)\right] \\
&=&S^{-1}(\Lambda (\gamma \rightharpoonup p_{L}^{2}))\cdot \left[
(a,b)(p_{L}^{1}\cdot c)\right] \\
&=&S^{-1}(\Lambda (\gamma \rightharpoonup p_{L}^{2}))\cdot \left[
\varepsilon (p_{L1}^{1})(a,b)(p_{L2}^{1}\cdot c)\right] \\
&=&S^{-1}(\Lambda (\gamma \rightharpoonup p_{L}^{2}))\cdot \left[
(p_{L1}^{1}\cdot (a,b))(p_{L2}^{1}\cdot c)\right] \\
&=&S^{-1}(\Lambda (\gamma \rightharpoonup p_{L}^{2}))p_{L}^{1}\cdot ((a,b)c)
\\
&=&S^{-1}(q_{L}^{1}p_{L1}^{2})p_{L}^{1}\cdot ((a,b)c)\gamma
(q_{L}^{2}p_{L2}^{2}) \\
&=&(a,b)c
\end{eqnarray*}%
for all $a,b,c\in A$, where in the fourth line we used the fact that $t$ is
a left integral and in the following lines that $(a,b)\in A^{H}$.
\end{itemize}
\end{proof}

\subsection{Surjectivity of the Morita maps}

\subsubsection{Galois Extensions}

\begin{definition}
Let $H$ be a finite dimensional quasi-Hopf algebra, $\dim _{k}H=n$, $A$ a
left $H$-module algebra and $B=A^{H}$ the subalgebra of $H$-invariants. We
say that the extension $B\subseteq A$ is \textbf{Galois }if the linear map $%
can:A\otimes _{B}A\longrightarrow A\otimes H^{\ast }$, $can(a\otimes
_{B}b)=\sum_{i=1}^{n}(p_{R}^{1}\cdot a)(p_{R}^{2}e_{i}\cdot b)\otimes e^{i}$
is bijective, where $(e_{i})_{i=1,n}$ and $(e^{i})_{i=1,n}$ are dual bases
for $H$, respectively $H^{\ast }$.
\end{definition}

\begin{remark}
As in the Hopf case, one may instead use the map $can^{\prime }:A\otimes
_{B}A\longrightarrow A\otimes H^{\ast }$, $can^{\prime }(a\otimes
_{B}b)=\sum_{i=1}^{n}(U_{R}^{1}e_{i}\cdot a)(U_{R}^{2}\cdot b)\otimes e^{i}$
where $U_{R}$ is the element from (\ref{UR}). If we define 
\begin{equation*}
\Xi :A\otimes H^{\ast }\longrightarrow A\otimes H^{\ast },\Xi (a\otimes
h^{\ast })=\sum_{i=1}^{n}e_{i}\cdot a\otimes (e^{i}\leftharpoonup
q_{L}^{1})(h^{\ast }S\leftharpoonup q_{L}^{2})
\end{equation*}%
then we have the connection between the two "$can$" maps: $\Xi \circ
can=can^{\prime }$. Also one may easily check that $\Xi $ is bijective, with
inverse given by 
\begin{equation*}
\Xi ^{-1}(a\otimes h^{\ast })=\sum_{i=1}^{n}e_{i}\cdot a\otimes
(f^{(-1)1}\rightharpoonup h^{\ast }S^{-1}\leftharpoonup
V_{L}^{1})(f^{(-1)2}\rightharpoonup e^{i}\leftharpoonup V_{L}^{2})
\end{equation*}%
where $f$ is the gauge element from (\ref{fdeltaf-1}) and $V_{L}$ is the
element introduced by (\ref{VL}).
\end{remark}

\begin{remark}
We proved in \textbf{Proposition} \ref{functorul indus} that the functor $%
(-)^{H}:$ $_{A\#H}\mathcal{M\longrightarrow }_{B}\mathcal{M}$ has a left
adjoint, namely $A\otimes _{B}(-)$. The counit of this adjunction is $%
\varepsilon _{M}:A\otimes _{B}M^{H}\longrightarrow M$, $\varepsilon
_{M}(a\otimes _{B}m)=am$ (for $M$ a left $A\#H$-module); it is left $A\#H$%
-linear. If $\varepsilon _{M}$ is an isomorphism for all modules $M$, we
call this the \textbf{Weak Structure Theorem}.
\end{remark}

\begin{remark}
Let $\theta _{t}:H^{\ast }\longrightarrow H$ be the isomorphism of the 
\textbf{Corollary }\ref{morfism frobenius}, $\theta _{t}(h^{\ast })=h^{\ast
}(q_{L}^{1}t_{1}p_{L}^{1})q_{L}^{2}t_{2}p_{L}^{2}$. As in the classical Hopf
case, we get the relation $(I_{A}\otimes \theta _{t})\circ can=\left[ -,-%
\right] $ which indicates that the Morita map $\left[ -,-\right] $ and the
Galois map $can$ will be simultaneously bijective:%
\begin{eqnarray*}
(I_{A}\otimes \theta _{t})\circ can(a\otimes _{B}b) &=&(I_{A}\otimes \theta
_{t})(\sum_{i=1}^{n}(p_{R}^{1}\cdot a)(p_{R}^{2}e_{i}\cdot b)\otimes e^{i})
\\
(\ref{morfism frobenius}) &=&\sum_{i=1}^{n}(p_{R}^{1}\cdot
a)(p_{R}^{2}e_{i}\cdot b)\otimes
e^{i}(q_{L}^{1}t_{1}p_{L}^{1})q_{L}^{2}t_{2}p_{L}^{2} \\
&=&(p_{R}^{1}\cdot a)(p_{R}^{2}q_{L}^{1}t_{1}p_{L}^{1}\cdot b)\otimes
q_{L}^{2}t_{2}p_{L}^{2} \\
(\ref{formpr}) &=&(x^{1}\cdot a)(x^{2}\beta
S(x^{3})q_{L}^{1}t_{1}p_{L}^{1}\cdot b)\otimes q_{L}^{2}t_{2}p_{L}^{2} \\
(\ref{Sqlt}) &=&(x^{1}\cdot a)(x^{2}\beta q_{L}^{1}t_{1}p_{L}^{1}\cdot
b)\otimes x^{3}q_{L}^{2}t_{2}p_{L}^{2} \\
(\ref{betaql}) &=&(x^{1}\cdot a)(x^{2}t_{1}p_{L}^{1}\cdot b)\otimes
x^{3}t_{2}p_{L}^{2} \\
(\ref{smashproduct}) &=&(a\#t)(p_{L}^{1}\cdot b\#p_{L}^{2}) \\
&=&\left[ a,b\right]
\end{eqnarray*}
\end{remark}

Now we have all the ingredients to prove the analogue of \textbf{Theorem}
3.1 of \cite{Beattie97a} in case of a finite dimensional quasi-Hopf algebra
and the proof follows closely the one in \cite{Beattie97a}.

\begin{theorem}
\label{ThMoritaGalois}Let $H$ be a finite dimensional quasi-Hopf algebra, $A$
a left $H$-module algebra and $B=A^{H}$ the subalgebra of $H$-invariants.
Then the following statements are equivalent:

\begin{enumerate}
\item The extension $A/B$ is Galois;

\item The map $can:A\otimes _{B}A\longrightarrow A\otimes H^{\ast }$ is
surjective;

\item The Morita map $\left[ -,-\right] $ is bijective;

\item The Morita map $\left[ -,-\right] $ is surjective;

\item The Weak Structure Theorem holds for $_{A}(_{H}\mathcal{M)}$;

\item The counit of the adjunction $\varepsilon _{M}:A\otimes
_{B}M^{H}\longrightarrow M$ is surjective for all left $A\#H$-modules $M$;

\item $A$ is a generator for the category $_{A}(_{H}\mathcal{M)\simeq }%
_{A\#H}\mathcal{M}$.
\end{enumerate}
\end{theorem}

\begin{proof}
1.$\Longrightarrow $2., 3.$\Longrightarrow $4., 5.$\Longrightarrow $6.
Obviously.

1.$\Longleftrightarrow $3., 2.$\Longleftrightarrow $4. Come from the
previous remark.

4.$\Longrightarrow $3. It results from the classical Morita theory for rings.

4.$\Longrightarrow $5. The injectivity of the counit of the adjunction: let $%
\sum_{i}a_{i}\otimes _{B}m_{i}\in A\otimes _{B}M^{H}$ such that $%
\sum_{i}a_{i}m_{i}=0$. By the surjectivity of $\left[ -,-\right] $ we may
find $\sum_{k}c_{k}\otimes _{B}d_{k}\in A\otimes _{B}A$ such that $\sum_{k}%
\left[ c_{k},d_{k}\right] =1_{A}\#1_{H}$ (which means $\sum_{k}(x^{1}\cdot
c_{k})(x^{2}t_{1}p_{L}^{1}\cdot d_{k})\#(x^{3}t_{2}p_{L}^{2})=1_{A}\#1_{H}$%
). We get then%
\begin{eqnarray}
\sum_{i}a_{i}\otimes _{B}m_{i} &=&\sum_{i}(1\#1)a_{i}\otimes _{B}m_{i} 
\notag \\
&=&\sum_{i,k}\left[ c_{k},d_{k}\right] a_{i}\otimes _{B}m_{i}  \notag \\
&=&\sum_{i,k}c_{k}(d_{k},a_{i})\otimes _{B}m_{i}  \notag \\
&=&\sum_{i,k}c_{k}\otimes _{B}(d_{k},a_{i})m_{i}  \notag \\
(\ref{moritaB}) &=&\sum_{i,k}c_{k}\otimes _{B}\left\{ t\left[
(p_{L}^{1}\cdot d_{k})(p_{L}^{2}\cdot a_{i})\right] \right\} m_{i}  \notag \\
&=&\sum_{i,k}c_{k}\otimes _{B}t\left\{ \left[ (p_{L}^{1}\cdot
d_{k})(p_{L}^{2}\cdot a_{i})\right] m_{i}\right\}   \notag \\
&=&\sum_{i,k}c_{k}\otimes _{B}t\left\{ (X^{1}p_{L}^{1}\cdot d_{k})\left[
(X^{2}p_{L}^{2}\cdot a_{i})(X^{3}m_{i})\right] \right\}   \notag \\
&=&\sum_{i,k}c_{k}\otimes _{B}t\left\{ (p_{L}^{1}\cdot d_{k})\left[
(p_{L}^{2}\cdot a_{i})m_{i}\right] \right\}   \notag \\
&=&\sum_{i,k}c_{k}\otimes _{B}t\left\{ (p_{L}^{1}\cdot d_{k})\left[
p_{L}^{2}\cdot (a_{i}m_{i})\right] \right\}   \notag \\
&=&0  \notag
\end{eqnarray}%
where in the last four lines we used that $m_{i}\in M^{H}$. For the
surjectivity, let $m\in M$. Then we compute that%
\begin{eqnarray*}
\varepsilon _{M}(\sum_{k}c_{k}\otimes _{B}t\left[ (p_{L}^{1}\cdot
d_{k})(p_{L}^{2}m)\right] ) &=&\sum_{k}c_{k}\left\{ t\left[ (p_{L}^{1}\cdot
d_{k})(p_{L}^{2}m)\right] )\right\}  \\
&=&\sum_{k}\underset{1_{A}\#1_{H}}{\underbrace{\left[ (x^{1}\cdot
c_{k})(x^{2}t_{1}p_{L}^{1}\cdot d_{k})\right] (x^{3}t_{2}p_{L}^{2}}}m) \\
&=&m
\end{eqnarray*}%
proving that $\varepsilon _{M}$ is indeed surjective.

5.$\Longrightarrow $3., 6.$\Longrightarrow $4. $A\#H$ is an $A\#H$-module
with the action given by the algebra multiplication, so by hypothesis $%
\varepsilon _{A\#H}$ is bijective (respectively surjective). We get the
following sequence of bijections (respectively surjections)%
\begin{equation*}
A\otimes _{B}A\simeq A\otimes _{B}(\int_{H}^{l}\otimes A)=A\otimes
_{B}(H\otimes A)^{H}\simeq A\otimes _{B}(A\#H)^{H}\overset{\varepsilon
_{A\#H}}{\longrightarrow }A\#H
\end{equation*}%
Explicitly, this means%
\begin{eqnarray*}
a\otimes _{B}b &\longrightarrow &a\otimes _{B}(t\otimes b)\longrightarrow
a\otimes _{B}(t_{1}p_{L}^{1}\cdot b\#t_{2}p_{L}^{2}) \\
&\longrightarrow &a\cdot (t_{1}p_{L}^{1}\cdot
b\#t_{2}p_{L}^{2})=(a\#1)(t_{1}p_{L}^{1}\cdot b\#t_{2}p_{L}^{2}) \\
&=&(x^{1}\cdot a)(x^{2}t_{1}p_{L}^{1}\cdot b)\#x^{3}t_{2}p_{L}^{2} \\
&=&\left[ a,b\right]
\end{eqnarray*}%
hence the Morita map $\left[ -,-\right] $ is bijective (respectively
surjective).

5.$\Longrightarrow $7, 7.$\Longrightarrow $6. The proofs are identical to
those in the Hopf case, so we omit them.
\end{proof}

We give now two examples of Galois extensions.

First, let $F\in H\otimes H$ be a gauge transformation. If we denote by $%
H_{F}$ the quasi-Hopf algebra obtained by twisting the comultiplication of $%
H $ via $F$, then in \cite{BulacuPanaite00} it is proven that there is a new
multiplication on $A$, namely $a\circ b=(F^{(-1)1}\cdot a)(F^{(-1)2}\cdot b)$%
, for $a,b\in A$ such that $A$, with this new multiplication (denoted $%
A_{F^{-1}}$), becomes a left $H_{F}$-module algebra. In this case the
categories $_{A}(_{H}\mathcal{M)}$ and $_{A_{F^{-1}}}(_{H_{F}}\mathcal{M)}$
are isomorphic and there is an algebra isomorphism between the smash
products $A\#H$ and $A_{F^{-1}}\#H_{F}$, which sends $a\#h\longrightarrow
F^{1}\cdot a\#F^{2}h$, for all $a\in A$, $h\in H$.

Remark also that $B$, the space of $H$-invariants, remains the same, as the
action of $H$ is not modified. Moreover, it is an associative algebra with
the multiplication induced by the new multiplication on $A_{F^{-1}}$, as $H$
acts trivially on $B$ and $F$ is a gauge transformation. As the categories $%
_{A}(_{H}\mathcal{M)}$ and $_{A_{F^{-1}}}(_{H_{F}}\mathcal{M)}$ are
isomorphic and the following diagram of functors is commutative%
\begin{equation*}
\begin{array}{ccccc}
_{A}(_{H}\mathcal{M)} &  & \simeq &  & _{A_{F^{-1}}}(_{H_{F}}\mathcal{M)} \\ 
& _{A\otimes _{B}(-)}\nwarrow \searrow ^{(-)^{H}} &  & ^{(-)^{H_{F}}}%
\nearrow \swarrow _{A_{F^{-1}}\otimes _{B}-} &  \\ 
&  & _{B}\mathcal{M} &  & 
\end{array}%
\end{equation*}%
as it can be easily checked, the counits of these two adjunctions will be
simultaneously bijective. Hence, by \textbf{Theorem} \ref{ThMoritaGalois}, $%
B\subseteq A$ is Galois $\Leftrightarrow $ $B\subseteq A_{F^{-1}}$ is Galois.

For the second example, let $\mathcal{A}$ be a right $H$-comodule algebra,
as it was defined in \cite{Hausser99}. That is, $\mathcal{A}$ is an
associative algebra endowed with an algebra morphism $\rho :\mathcal{%
A\longrightarrow A\otimes }H$ and an invertible element $\phi _{\rho }\in 
\mathcal{A\otimes }H\mathcal{\otimes }H$, such that 
\begin{subequations}
\begin{eqnarray}
\phi _{\rho }(\rho \otimes I)\rho (a)\phi _{\rho }^{-1} &=&(I\otimes \Delta
)\rho (a)  \label{coasoccomod} \\
(I\otimes \varepsilon )\rho (a) &=&a  \label{counitcomod} \\
(I\otimes I\otimes \Delta )(\phi _{\rho })(\rho \otimes I\otimes I)(\phi
_{\rho }) &=&(1_{\mathcal{A}}\otimes \phi )(I\otimes \Delta \otimes I)(\phi
_{\rho })(\phi _{\rho }\otimes 1)  \label{pentagoncomod} \\
(I\otimes \varepsilon \otimes I)(\phi _{\rho }) &=&(I\otimes I\otimes
\varepsilon )(\phi _{\rho })=1\otimes 1  \label{counitrocomod}
\end{eqnarray}%
hold for all $a\in \mathcal{A}$. We shall denote $\phi _{\rho }=X_{\rho
}^{1}\otimes X_{\rho }^{2}\otimes X_{\rho }^{3}$ and $\phi _{\rho
}^{-1}=x_{\rho }^{1}\otimes x_{\rho }^{2}\otimes x_{\rho }^{3}$. Following 
\cite{Bulacu02}, we may define the quasi-smash product $\mathcal{A}\overline{%
\#}H^{\ast }$. As vector space, this is $\mathcal{A}\otimes H^{\ast }$
endowed with a multiplication given by 
\end{subequations}
\begin{equation}
(a\overline{\#}h^{\ast })(b\overline{\#}g^{\ast })=ab_{0}x_{\rho }^{1}%
\overline{\#}(h^{\ast }\leftharpoonup b_{1}x_{\rho }^{2})(g^{\ast
}\leftharpoonup x_{\rho }^{3})  \label{quasismash}
\end{equation}%
for any $a$, $b\in \mathcal{A}$, $h^{\ast }$, $g^{\ast }\in H^{\ast }$.
Using the left $H$-action given by 
\begin{equation*}
g(a\overline{\#}h^{\ast })=a\overline{\#}g\rightharpoonup h^{\ast }
\end{equation*}%
for any $g\in H$, $a\in \mathcal{A}$, $h^{\ast }\in H^{\ast }$, the
quasi-smash product $\mathcal{A}\overline{\#}H^{\ast }$ becomes a left $H$%
-module algebra with invariants $(\mathcal{A}\overline{\#}H^{\ast })^{H}=%
\mathcal{A}\overline{\#}k\varepsilon \simeq \mathcal{A}$. For this module
algebra, the Galois map $can$ is 
\begin{eqnarray*}
can &:&(\mathcal{A}\overline{\#}H^{\ast })\otimes _{\mathcal{A}}(\mathcal{A}%
\overline{\#}H^{\ast })\longrightarrow (\mathcal{A}\overline{\#}H^{\ast
})\otimes H^{\ast } \\
(a\overline{\#}h^{\ast })\otimes _{\mathcal{A}}(b\overline{\#}g^{\ast })
&\longrightarrow &\sum_{i=1}^{n}(a\overline{\#}p_{R}^{1}\rightharpoonup
h^{\ast })(b\overline{\#}p_{R}^{2}e_{i}\rightharpoonup g^{\ast })\otimes
e^{i}
\end{eqnarray*}%
where $a$, $b\in \mathcal{A}$, $h^{\ast }$, $g^{\ast }\in H^{\ast }$. But $%
\mathcal{A}\overline{\#}H^{\ast }=(\mathcal{A}\overline{\#}\varepsilon )(1_{%
\mathcal{A}}\overline{\#}H^{\ast })$, meaning that it's enough to consider
elements of the type $(a\overline{\#}h^{\ast })\otimes _{\mathcal{A}}(1_{%
\mathcal{A}}\overline{\#}g^{\ast })$, as in the Hopf case. Now the formula
for the Galois map becomes 
\begin{equation*}
can((a\overline{\#}h^{\ast })\otimes _{\mathcal{A}}(1_{\mathcal{A}}\overline{%
\#}g^{\ast }))=\sum_{i=1}^{n}(a\overline{\#}p_{R}^{1}\rightharpoonup h^{\ast
})(1_{\mathcal{A}}\overline{\#}p_{R}^{2}e_{i}\rightharpoonup g^{\ast
})\otimes e^{i}
\end{equation*}%
We need the following element introduced in \cite{Hausser99}: $q_{\rho
}=q_{\rho }^{1}\otimes q_{\rho }^{2}=X_{\rho }^{1}\otimes S^{-1}(\alpha
X_{\rho }^{3})X_{\rho }^{2}$. This element has similar properties with $%
q_{R} $:%
\begin{eqnarray}
(1_{\mathcal{A}}\otimes S^{-1}(a_{1}))q_{\rho }\rho (a_{0}) &=&(a\otimes
1)q_{\rho }  \label{qro} \\
(\rho \otimes I)(q_{R})\phi _{\rho }^{-1} &=&(1_{\mathcal{A}}\otimes
1\otimes S^{-1}(X_{\rho }^{3}))(X_{\rho }^{1}\otimes q_{\rho }\Delta
(X_{\rho }^{2}))  \label{qrofi}
\end{eqnarray}%
for all $a\in \mathcal{A}$.

\begin{proposition}
The quasi-smash product $\mathcal{A}\overline{\#}H^{\ast }$ is a Galois
extension of $\mathcal{A}$, with inverse of the Galois map given by%
\begin{equation*}
can^{-1}((a\overline{\#}h^{\ast })\otimes g^{\ast })=\sum_{i=1}^{n}(a%
\overline{\#}h^{\ast })(q_{\rho }^{1}\overline{\#}e^{i}S\leftharpoonup
q_{\rho }^{2})\otimes _{\mathcal{A}}(1_{\mathcal{A}}\overline{\#}g^{\ast
}\leftharpoonup e_{i})
\end{equation*}%
for all $a\in \mathcal{A}$, $h^{\ast }$, $g^{\ast }\in H^{\ast }$.
\end{proposition}

\begin{proof}
For all $a\in \mathcal{A}$ and $h^{\ast }$, $g^{\ast }\in H^{\ast }$, we
compute that%
%TCIMACRO{\TeXButton{breakdisplay}{\begin{allowdisplaybreaks}}}%
%BeginExpansion
\begin{allowdisplaybreaks}%
%EndExpansion
\begin{eqnarray*}
can^{-1}\circ can((a\overline{\#}h^{\ast })\otimes _{\mathcal{A}}(1_{%
\mathcal{A}}\overline{\#}g^{\ast })) &=&can^{-1}(\sum_{i=1}^{n}(a\overline{\#%
}p_{R}^{1}\rightharpoonup h^{\ast })(1_{\mathcal{A}}\overline{\#}%
p_{R}^{2}e_{i}\rightharpoonup g^{\ast })\otimes e^{i}) \\
&=&\sum_{i,j=1}^{n}[(a\overline{\#}p_{R}^{1}\rightharpoonup h^{\ast })(1_{%
\mathcal{A}}\overline{\#}p_{R}^{2}e_{i}\rightharpoonup g^{\ast })](q_{\rho
}^{1}\overline{\#}e^{j}S\leftharpoonup q_{\rho }^{2})\otimes _{\mathcal{A}%
}(1_{\mathcal{A}}\overline{\#}e^{i}\leftharpoonup e_{j}) \\
&=&\sum_{i,j=1}^{n}(a\overline{\#}X^{1}p_{R}^{1}\rightharpoonup h^{\ast
})[(1_{\mathcal{A}}\overline{\#}X^{2}p_{R}^{2}e_{j}e_{i}\rightharpoonup
g^{\ast })(q_{\rho }^{1}\overline{\#}X^{3}\rightharpoonup
e^{j}S\leftharpoonup q_{\rho }^{2})] \\
&&\otimes _{\mathcal{A}}(1_{\mathcal{A}}\overline{\#}e^{i}) \\
&=&\sum_{i,j=1}^{n}(a\overline{\#}X^{1}p_{R}^{1}\rightharpoonup h^{\ast
})[(1_{\mathcal{A}}\overline{\#}X^{2}p_{R}^{2}S(X^{3})e_{j}S(q_{\rho
}^{2})e_{i}\rightharpoonup g^{\ast })(q_{\rho }^{1}\overline{\#}e^{j}S)] \\
&&\otimes _{\mathcal{A}}(1_{\mathcal{A}}\overline{\#}e^{i}) \\
(\ref{formpr}) &=&\sum_{i,j=1}^{n}(a\overline{\#}h^{\ast })[(1_{\mathcal{A}}%
\overline{\#}\beta e_{j}S(q_{\rho }^{2})e_{i}\rightharpoonup g^{\ast
})(q_{\rho }^{1}\overline{\#}e^{j}S)]\otimes _{\mathcal{A}}(1_{\mathcal{A}}%
\overline{\#}e^{i}) \\
(\ref{quasismash}) &=&\sum_{i,j=1}^{n}(a\overline{\#}h^{\ast })(q_{\rho
_{0}}^{1}x_{\rho }^{1}\overline{\#}(\beta e_{j}S(q_{\rho
}^{2})e_{i}\rightharpoonup g^{\ast }\leftharpoonup q_{\rho _{1}}^{1}x_{\rho
}^{2})(e^{j}S\leftharpoonup x_{\rho }^{3})) \\
&&\otimes _{\mathcal{A}}(1_{\mathcal{A}}\overline{\#}e^{i}) \\
&=&\sum_{i,j=1}^{n}(a\overline{\#}h^{\ast })(q_{\rho _{0}}^{1}x_{\rho }^{1}%
\overline{\#}(\beta e_{j}S(x_{\rho }^{3})S(q_{\rho
}^{2})e_{i}\rightharpoonup g^{\ast }\leftharpoonup q_{\rho _{1}}^{1}x_{\rho
}^{2})(e^{j}S)) \\
&&\otimes _{\mathcal{A}}(1_{\mathcal{A}}\overline{\#}e^{i}) \\
(\ref{pentagoncomod}) &=&\sum_{i,j=1}^{n}(a\overline{\#}h^{\ast })(x_{\rho
}^{1}X_{\rho }^{1}\overline{\#}(\beta e_{j}S(x_{\rho _{1}}^{3}X^{2}X_{\rho
_{2}}^{2})\alpha x_{\rho _{2}}^{3}X^{3}X_{\rho }^{3}e_{i}\rightharpoonup
g^{\ast } \\
&\leftharpoonup &x_{\rho }^{2}X^{1}X_{\rho _{1}}^{2})(e^{j}S))\otimes _{%
\mathcal{A}}(1_{\mathcal{A}}\overline{\#}e^{i}) \\
(\ref{alfa}),(\ref{counitrocomod}) &=&\sum_{i,j=1}^{n}(a\overline{\#}h^{\ast
})(X_{\rho }^{1}\overline{\#}(\beta e_{j}S(X^{2}X_{\rho _{2}}^{2})\alpha
X^{3}X_{\rho }^{3}e_{i}\rightharpoonup g^{\ast } \\
&\leftharpoonup &X^{1}X_{\rho _{1}}^{2})(e^{j}S))\otimes _{\mathcal{A}}(1_{%
\mathcal{A}}\overline{\#}e^{i}) \\
(\ref{beta}) &=&\sum_{i=1}^{n}(a\overline{\#}h^{\ast })(X_{\rho }^{1}%
\overline{\#}g^{\ast }(X^{1}X_{\rho _{1}}^{2}\beta S(X^{2}X_{\rho
_{2}}^{2})\alpha X^{3}X_{\rho }^{3}e_{i})\varepsilon ) \\
&&\otimes _{\mathcal{A}}(1_{\mathcal{A}}\overline{\#}e^{i}) \\
(\ref{beta}) &=&\sum_{i=1}^{n}(a\overline{\#}h^{\ast })(1_{\mathcal{A}}%
\overline{\#}g^{\ast }(X^{1}\beta S(X^{2})\alpha X^{3}e_{i})\varepsilon
)\otimes _{\mathcal{A}}(1_{\mathcal{A}}\overline{\#}e^{i}) \\
(\ref{fibetaalfa}) &=&\sum_{i=1}^{n}(a\overline{\#}h^{\ast })(1_{\mathcal{A}}%
\overline{\#}g^{\ast }(e_{i})\varepsilon )\otimes _{\mathcal{A}}(1_{\mathcal{%
A}}\overline{\#}e^{i}) \\
&=&\sum_{i=1}^{n}(a\overline{\#}h^{\ast })\otimes _{\mathcal{A}}(1_{\mathcal{%
A}}\overline{\#}g^{\ast })
\end{eqnarray*}%
\pagebreak 
%TCIMACRO{\TeXButton{endbreakdisplay}{\end{allowdisplaybreaks}}}%
%BeginExpansion
\end{allowdisplaybreaks}%
%EndExpansion
Next, we have that

%TCIMACRO{\TeXButton{breakdisplay}{\begin{allowdisplaybreaks}}}%
%BeginExpansion
\begin{allowdisplaybreaks}%
%EndExpansion
\begin{eqnarray*}
can\circ can^{-1}((a\overline{\#}h^{\ast })\otimes g^{\ast })
&=&can(\sum_{j=1}^{n}(a\overline{\#}h^{\ast })(q_{\rho }^{1}\overline{\#}%
e^{j}S\rightharpoonup q_{\rho }^{2})\otimes _{\mathcal{A}}(1_{\mathcal{A}}%
\overline{\#}g^{\ast }\leftharpoonup e_{j})) \\
(\ref{quasismash}) &=&can(\sum_{j=1}^{n}(aq_{\rho _{0}}^{1}x_{\rho }^{1}%
\overline{\#}(h^{\ast }\leftharpoonup q_{\rho _{1}}^{1}x_{\rho
}^{2})(e^{j}S\rightharpoonup q_{\rho }^{2}x_{\rho }^{3}))\otimes _{\mathcal{A%
}}(1_{\mathcal{A}}\overline{\#}g^{\ast }\leftharpoonup e_{j})) \\
&=&\sum_{i,j=1}^{n}[aq_{\rho _{0}}^{1}x_{\rho }^{1}\overline{\#}%
(p_{R_{1}}^{1}\rightharpoonup h^{\ast }\leftharpoonup q_{\rho
_{1}}^{1}x_{\rho }^{2})(p_{R_{2}}^{1}\rightharpoonup e^{j}S\rightharpoonup
q_{\rho }^{2}x_{\rho }^{3})](1_{\mathcal{A}}\overline{\#}p_{R}^{2}e_{i}%
\rightharpoonup g^{\ast } \\
&\leftharpoonup &e_{j})\otimes e^{i} \\
(\ref{quasismash}) &=&\sum_{i,j=1}^{n}(aq_{\rho _{0}}^{1}x_{\rho
}^{1}y_{\rho }^{1}\overline{\#}[(p_{R_{1}}^{1}\rightharpoonup h^{\ast
}\leftharpoonup q_{\rho _{1}}^{1}x_{\rho }^{2}y_{\rho
_{1}}^{2})(p_{R_{2}}^{1}\rightharpoonup e^{j}S\rightharpoonup q_{\rho
}^{2}x_{\rho }^{3}y_{\rho _{2}}^{2})](p_{R}^{2}e_{i}\rightharpoonup g^{\ast }
\\
&\leftharpoonup &e_{j}y_{\rho }^{3}))\otimes e^{i} \\
(\ref{pentagoncomod}) &=&\sum_{i,j=1}^{n}(ax_{\rho }^{1}\overline{\#}%
(X^{1}p_{R_{1}}^{1}\rightharpoonup h^{\ast }\leftharpoonup x_{\rho
}^{2})[(X^{2}p_{R_{2}}^{1}\rightharpoonup e^{j}S\rightharpoonup
S^{-1}(\alpha )x_{\rho _{1}}^{3})(X^{3}p_{R}^{2}e_{i}\rightharpoonup g^{\ast
} \\
&\leftharpoonup &e_{j}x_{\rho _{2}}^{3})])\otimes e^{i} \\
(\ref{alfa}) &=&\sum_{i,j=1}^{n}(a\overline{\#}(X^{1}p_{R_{1}}^{1}%
\rightharpoonup h^{\ast })[(X^{2}p_{R_{2}}^{1}\rightharpoonup
e^{j}S)(X^{3}p_{R}^{2}e_{i}\rightharpoonup g^{\ast }\leftharpoonup
e_{j}\alpha )])\otimes e^{i} \\
&=&\sum_{i=1}^{n}a\overline{\#}(X^{1}p_{R_{1}}^{1}\rightharpoonup h^{\ast
})\otimes g^{\ast }(S(X^{2}p_{R_{2}}^{1})\alpha X^{3}p_{R}^{2}e_{i})e^{i} \\
(\ref{prqr}) &=&a\overline{\#}h^{\ast }\otimes g^{\ast }
\end{eqnarray*}%
%TCIMACRO{\TeXButton{endbreakdisplay}{\end{allowdisplaybreaks}}}%
%BeginExpansion
\end{allowdisplaybreaks}%
%EndExpansion
\end{proof}

\subsubsection{Total integrals}

\begin{definition}
Let $H$ be a quasi-Hopf algebra and $A$ a left $H$-module algebra. A \textbf{%
total integral} for $A$ is a left $H$-morphism $\Phi :H^{\ast
}\longrightarrow A$ such that $\Phi (\varepsilon )=1_{A}$ (on $H^{\ast }$ we
take the structure of left $H$-module given by translation: $%
(h\rightharpoonup h^{\ast })(g)=h^{\ast }(gh)$ for all $h^{\ast }\in H^{\ast
}$, $g,h\in H$).
\end{definition}

\begin{proposition}
With notations as above, the following statements are equivalent:

\begin{enumerate}
\item The Morita map $\left( -,-\right) $ is surjective;

\item There is a total integral for $A$;

\item $A$ has an element of trace one (i.e. $a\in A$ such that $t\cdot
a=1_{A}$).
\end{enumerate}
\end{proposition}

\begin{proof}
1.$\Longrightarrow $2. Let $\sum_{i}a_{i}\otimes _{B}b_{i}\in A\otimes _{B}A$
such that $\sum_{i}(a_{i},b_{i})=1_{A}$. We define $\Phi :H^{\ast
}\longrightarrow A$, 
\begin{equation*}
\Phi (h^{\ast })=\sum\nolimits_{i}q_{L}^{1}t_{1}\cdot \left[ (p_{L}^{1}\cdot
a_{i})(p_{L}^{2}\cdot b_{i})\right] h^{\ast }S(q_{L}^{2}t_{2})
\end{equation*}%
Then 
\begin{eqnarray*}
\Phi (\varepsilon ) &=&\sum\nolimits_{i}q_{L}^{1}t_{1}\cdot \left[
(p_{L}^{1}\cdot a_{i})(p_{L}^{2}\cdot b_{i})\right] \varepsilon
(q_{L}^{2}t_{2}) \\
(\ref{formql}) &=&\sum\nolimits_{i}\alpha t\cdot \left[ (p_{L}^{1}\cdot
a_{i})(p_{L}^{2}\cdot b_{i})\right] \\
&=&t\cdot \left[ (p_{L}^{1}\cdot a_{i})(p_{L}^{2}\cdot b_{i})\right] \\
&=&\sum_{i}(a_{i},b_{i})=1_{A}
\end{eqnarray*}%
and 
\begin{eqnarray*}
\Phi (h\rightharpoonup h^{\ast }) &=&\sum\nolimits_{i}q_{L}^{1}t_{1}\cdot 
\left[ (p_{L}^{1}\cdot a_{i})(p_{L}^{2}\cdot b_{i})\right] h^{\ast
}S(S^{-1}(h)q_{L}^{2}t_{2}) \\
(\ref{Sqlt}) &=&\sum\nolimits_{i}hq_{L}^{1}t_{1}\cdot \left[ (p_{L}^{1}\cdot
a_{i})(p_{L}^{2}\cdot b_{i})\right] h^{\ast }S(q_{L}^{2}t_{2}) \\
&=&h\cdot \Phi (h^{\ast })
\end{eqnarray*}

2.$\Longrightarrow $1. Let $\Phi $ be a total integral. As $t\neq 0$ and $H$
is finite dimensional, we may find $h^{\ast }\in H^{\ast }$ such that $%
h^{\ast }(t)=1$. Then $(t\rightharpoonup h^{\ast })(h)=h^{\ast
}(ht)=\varepsilon (h)h^{\ast }(t)=\varepsilon (h)$, $\forall h\in H$,
meaning that $t\rightharpoonup h^{\ast }=\varepsilon $. Hence%
\begin{eqnarray*}
\left( 1_{A},\Phi (h^{\ast })\right) &=&t\cdot \left[ (p_{L}^{1}\cdot
1_{A})(p_{L}^{2}\cdot \Phi (h^{\ast }))\right] \\
&=&t\cdot \Phi (h^{\ast }) \\
&=&\Phi (t\rightharpoonup h^{\ast }) \\
&=&\Phi (\varepsilon )=1_{A}
\end{eqnarray*}%
and using the $B$-bilinearity of $\left( -,-\right) $ we get the
surjectivity.

2.$\Longrightarrow $3. Let $\Phi $ be a total integral. As above, consider $%
h^{\ast }\in H^{\ast }$ such that $t\rightharpoonup h^{\ast }=\varepsilon $.
Then $a=t\cdot \Phi (h^{\ast })$ is a trace one element.

3.$\Longrightarrow $2. Let $a\in A$ an element of trace one and $h^{\ast
}\in H^{\ast }$, $t\rightharpoonup h^{\ast }=\varepsilon $ as above. Define $%
\Phi :H^{\ast }\longrightarrow A$, $\Phi (g^{\ast })=q_{L}^{1}t_{1}\cdot
ag^{\ast }S(q_{L}^{2}t_{2})$. Then $\Phi (\varepsilon )=q_{L}^{1}t_{1}\cdot
a\varepsilon S(q_{L}^{2}t_{2})=\alpha t\cdot a=t\cdot a=1_{A}$ and $\Phi
(h\rightharpoonup g^{\ast })=q_{L}^{1}t_{1}\cdot ag^{\ast
}S(S^{-1}(h)q_{L}^{2}t_{2})=hq_{L}^{1}t_{1}\cdot ag^{\ast
}S(q_{L}^{2}t_{2})=h\Phi (g^{\ast })$, which means that $\Phi $ is a total
integral.
\end{proof}

\begin{example}
1) Let $F\in H\otimes H$ be a gauge transformation, as in the example of the
previous section. Then a total integral for $A$ remains a total integral for 
$A_{F^{-1}}$, as the action of $H$ is not modified. So the Morita maps $%
(-,-) $ will be simultaneously bijective.

2) Let $\mathcal{A}$ be a right $H$-comodule algebra. Then it is easy to see
that the map $\Phi :H^{\ast }\longrightarrow \mathcal{A}\overline{\#}H^{\ast
}$, given by $\Phi (h^{\ast })=1\overline{\#}h^{\ast }$, is a total
integral. Hence in this example we get the equivalence of categories $_{%
\mathcal{A}}\mathcal{M\simeq }_{(\mathcal{A}\overline{\#}H^{\ast })\#H}%
\mathcal{M}$
\end{example}

In connection with the notion of total integral, Bulacu and Nauwelaerts
proved in \cite{Bulacu00} the following three statements, for a dual
quasi-Hopf algebra and a comodule algebra. But in the finite dimensional
case this is the same as working with the quasi-Hopf algebra and a module
algebra, so we state them for completeness:

\begin{proposition}
(Proposition 2.9, \cite{Bulacu00}) If $A$ is a left $H$-module algebra and
there is a total integral $\Phi :H^{\ast }\longrightarrow A$, then each
relative module $M\in (_{H}\mathcal{M)}_{A}$ is injective as an $H$-module
(where $(_{H}\mathcal{M)}_{A}$ is the category of right $A$-modules in the
monoidal category $_{H}\mathcal{M}$).
\end{proposition}

\begin{corollary}
(Corollary 2.10, \cite{Bulacu00}) Under the previous hypotheses, the
following statements are equivalent:

\begin{enumerate}
\item $A$ is an injective left $H$-module;

\item There is a total integral on $A$;

\item Each object in $(_{H}\mathcal{M)}_{A}$ is injective as an $H$-module.
\end{enumerate}
\end{corollary}

\begin{theorem}
(Theorem 2.11, \cite{Bulacu00}) If $A$ is a left $H$-module algebra and
there is a total integral $\Phi :H^{\ast }\longrightarrow A$ which is
multiplicative, then for every $M\in (_{H}\mathcal{M)}_{A}$ the counit $%
\varepsilon _{M}:A\otimes _{B}M^{H}\longrightarrow M$, $\varepsilon
_{M}(a\otimes _{B}m)=am$, of the adjunction $\mathcal{M}_{B}\overset{%
\mathcal{-}\otimes _{B}A}{\underset{(-)^{H}}{\mathcal{\rightleftarrows }}}$ $%
(_{H}\mathcal{M)}_{A}$ is an isomorphism.
\end{theorem}

As remarked in the quoted paper, working with left or right $A$-modules is
essentially the same, just by passing to the opposite algebra $A^{op}$
(which is a module algebra over $H^{op,cop}$), so we can rephrase these
results in our context and obtain the following theorem:

\begin{theorem}
Let $H$ be a finite dimensional quasi-Hopf algebra, $A$ a left $H$-module
algebra and $B=A^{H}$ the subalgebra of $H$-invariants. Then the following
statements are equivalent:

\begin{enumerate}
\item The Morita map $\left( -,-\right) $ is surjective;

\item There is a total integral for $A$;

\item $A$ has an element of trace one (i.e. $a\in A$ such that $t\cdot
a=1_{A}$);

\item $A$ is an injective left $H$-module;

\item Each object in $_{A}(_{H}\mathcal{M)}$ is injective as an $H$-module.
\end{enumerate}
\end{theorem}

Combining the results of the previous two sections, we may state now the
following:

\begin{theorem}
Let $H$ be a finite dimensional quasi-Hopf algebra, $A$ a left $H$-module
algebra and $B=A^{H}$ the subalgebra of $H$-invariants. Then the following
statements are equivalent:

\begin{enumerate}
\item The functors $_{B}\mathcal{M}\overset{A\otimes _{B}(-)}{\underset{%
(-)^{H}}{\mathcal{\rightleftarrows }}}$ $_{A}(_{H}\mathcal{M)}$ are a pair
of inverse equivalences (\textbf{Strong Structure Theorem)};

\item The Morita maps $[,]$ and $(,)$ are surjective;

\item The extension $B\subseteq A$ is Galois and there is a total integral
on $A$.
\end{enumerate}
\end{theorem}

\begin{conclusion}
As noticed in the introduction, if we want to start with a right $H$%
-comodule algebra $\mathcal{A}$ and take the usual definition for
coinvariants $\mathcal{A}^{coH}$, this does not work any more in the
quasi-Hopf setting. For example, if we take a left $H$-module algebra $A$
and form the smash product $A\#H$, then this is a right comodule algebra (%
\cite{BulacuPanaite00}), but we cannot recover $A$ from the coinvariants as
in the Hopf case (we get something bigger). There are two possible ways to
overcome this problem: either to find an adequate definition for the
coinvariants, as it was done in \cite{Panaite05}, or to pass to
bicategories. Anyway, if we want the previous example to fit, we need some
coinvariants which are associative in the category of left $H$-modules, so
we can only get a Morita context in this monoidal category. It would be
interesting to see which are the connections between these two types of
Morita contexts, knowing that in the Hopf case these two contexts are the
same.
\end{conclusion}

\begin{center}
\textbf{ACKNOWLEDGMENT}
\end{center}

The author would like to thank Prof. C. N\u{a}st\u{a}sescu, F. Panaite and
D. Bulacu for their useful comments which improved this paper.

\bibliographystyle{plain}
\bibliography{ADRIANA}

Department of Mathematics, University Politehnica of Bucharest, 313 Splaiul
Independen\c{t}ei, 060042 Bucharest, Romania

\textit{E-mail address:} steleanu@mailcity.com

\end{document}